\documentclass{amsart}
\usepackage[authoryear]{natbib}



\begin{document}

\title[sequentially ordered hidden
Markov models]
{On the equivalence between standard and sequentially ordered hidden
Markov models
\footnote{This paper has been published in Stats and Prob. Letters (2008, vol. 78, pages 2171-2174).}}

\author{N. CHOPIN}
\address{CREST-ENSAE, Timbre J120 - 3, Avenue Pierre Larousse - 92245 Malakoff CEDEX,
FRANCE}
\email{nicolas.chopin@ensae.fr}

\maketitle
\begin{abstract}
\citet{Chopin2007} introduced a sequentially ordered hidden Markov
model, for which states are ordered according to their order of appearance,
and claimed that such a model is a re-parametrisation of a standard
Markov model. This note gives a formal proof that this equivalence
holds in Bayesian terms, as both formulations generate equivalent
posterior distributions, but does not hold in Frequentist terms, as
both formulations generate incompatible likelihood functions. Perhaps
surprisingly, this shows that Bayesian re-parametrisation and Frequentist
re-parametrisation are not identical concepts.

\end{abstract}
A standard hidden  Markov model assumes that the observed process
${(y}_{t})$ follows a mixture distribution:\begin{equation}
y_{t}|\{s_{t}=k\}\sim f(y_{t}|\xi_{k})\qquad t=0,1,\ldots\label{eq:obs}\end{equation}
where $\left\{ f(\cdot|\xi)\,\xi\in\Xi\right\} $ is a parametric
family of probability densities, and that the unobserved assignment
variable $(s_{t})$ is a $K$-state Markov chain: \begin{equation}
p(s_{t+1}=l|s_{t}=k)=q_{kl},\qquad k,l=1,\ldots,K.\label{eq:stdsys}\end{equation}

If necessary, distribution (\ref{eq:obs}) can also depends on $y_{t-1},y_{t-2},\ldots$
The transition matrix is denoted by $Q=(q_{kl})$, and the vector
of unknown parameters, which includes both the $q_{kl}$'s and the
$\xi_{k}$'s, is denoted by $\theta$. Lastly, it is often assumed
that\begin{equation}
p(s_{0}=k)=\vartheta_{k},\qquad k=1,\ldots,K,\label{eq:s0}\end{equation}
where $\vartheta=(\vartheta_{1},\ldots,\vartheta_{K})$ is the the
stationary distribution of the chain, i.e. the solution of $\vartheta Q=\vartheta$.
See e.g. \citet{McZuc}, \citet[Chap. 13]{MclachlanPeelMixtures},
\citet{Scott:HMM} or \citet{CapMouRyd} for more background on hidden
Markov models and their numerous applications. 

The sequentially ordered hidden Markov model introduced by \citet{Chopin2007}
has the same observation equation as (\ref{eq:obs}), i.e. 

\begin{equation}
y_{t}|\{z_{t}=k\}\sim f(y_{t}|\overline{\xi}_{k})\qquad t=0,1,\ldots\label{eq:obs2}\end{equation}

but differs with respect to the behaviour of the hidden process, now
denoted by $(z_{t},m_{t})$: \[
(z_{0},m_{0})=(1,1),\]
\begin{equation}
p(z_{t+1}=l|z_{t}=k,m_{t}=m)=\left\{ \begin{array}{ll}
\overline{q}{}_{kl} & \mbox{if }k,l\leq m\leq K,\\
\sum_{i=m+1}^{K}\overline{q}_{ki} & \mbox{if }k<l=m+1\leq K,\\
0 & \mbox{otherwise,}\end{array}\right.\label{eq:sys2}\end{equation}

\[
m_{t+1}=\max(m_{t},z_{t+1}),\]

for $t=0,1,\ldots$ For sake of clarity, parameters in this second
formulation are barred, e.g. $\overline{\theta}$, $\overline{\xi}_{k}$,
$\overline{q}_{kl}$, etc. The quantity $m_{t}$ represents the number
of states that have appeared up to time $t$, and $z_{t}$ represents
the current state, as labelled with respect to the order of appearance
of state values; e.g. $z_{0}=1$, as first state to appear is always
labelled `1', then $z_{1}=1$ with probability $q_{11}$, and $z_{1}=2$
otherwise, etc. The pair $(z_{t},m_{t})$ is a $K'$-state Markov
chain, with $K'=K(K+1)/2$, since $z_{t}\leq m_{t}$ with probability
one. 

\citet{Chopin2007} discusses several advantages of time-ordered hidden
Markov models. First, they are identifiable, provided $T\geq K$,
whereas standard hidden Markov models are invariant with respect to
state re-labelling. Second, they are still hidden Markov models, with
hidden chain $(z_{t},m_{t})$, so standard algorithms for hidden Markov
models (such as Gibbs sampling) can be adapted with little extra effort.
Third, one may conveniently estimate $m_{t}$ in order to evaluate
how many states are required to model the data. Fourth, in sequential
settings, that is, where statistical inference is performed at each
time $t$ where a new data-point is available, states are automatically
and consistently ordered at all iterations. 

A last advantage of these models, which was not mentioned by \citet{Chopin2007}
is that they are slightly more parsimonious, since they do not require
the specification of a distribution for the initial state $s_{0}$,
as in (\ref{eq:s0}). 

The rest of the paper is organised as follows. Section 1 proves that
the two formulations generate equivalent posterior distributions.
Section 2 shows that the two formulations lead to likelihood functions
that cannot be compared with each other. Section 3 discusses this
paradox, and mentions possible extensions of sequentially ordered
hidden Markov models.

\section{Bayesian equivalence}

We assume that $(y_{t})$ is the stochastic process defined by (\ref{eq:obs})
and (\ref{eq:stdsys}), i.e. a standard hidden Markov model, where
dependencies on parameter $\theta$ are interpreted as conditional
dependencies on random variable $\theta$, with prior probability
density $\pi$. We prove that, under the following two mild conditions,
one may define a latent process $(z_{t},m_{t})$ and a random parameter
$\overline{\theta}$ distributed according to $\pi$, such that the
distribution of $(y_{t})$ conditional on $\overline{\theta}$ corresponds
to (\ref{eq:obs2}) and (\ref{eq:sys2}), i.e. a sequentially ordered
hidden Markov model. 

\textbf{Condition 1.} The prior distribution $\pi$ is invariant with
respect to state re-labelling, that is, for any permutation $\tau$
of the first $K$ integers, if $\theta\sim\pi$, then \[
\theta_{\tau}=(\xi_{\tau(1)},\ldots,\xi_{\tau(K)},q_{\tau(1)\tau(1)},\ldots,q_{\tau(K)\tau(K-1)})\]
 is also distributed according to $\pi$.


\textbf{Condition 2.} Under the prior distribution $\pi$, all the
components of matrix $Q$ are positive with probability one. 

Let $\sigma_{t}$ be the random sequence of increasing size that records
sequentially the state values as they appear for the first time in
$s_{1:t}$; e.g. if $x_{1:5}=(4,3,4,7,3)$ then $\sigma_{9}=(4,3,7)$.
With an abuse of notations, $\sigma_{t}$ shall also stand for any
permutation $\tau$ of the first $K$ integers such that $\tau(i)=\sigma_{t}(i)$,
for $1\leq i\leq k$, where $k$ is the length of sequence $\sigma_{t}$.
In particular, let $z_{t}=\sigma_{t}^{-1}(s_{t})$ and $m_{t}=\max_{1\leq t'\leq t}z_{t'}$

According to Condition 2 and basic properties of Markov chains, there
exists almost surely a vector $\sigma$ of size $K$ and a finite
time $\zeta$ such that $\sigma_{t}=\sigma$ for all $t\geq\zeta$.
Let $\overline{\theta}=\sigma(\theta)$ and $\overline{Q}=\sigma(Q)$,
which are obtained respectively by re-ordering the components of $\theta$
and $Q$ according to $\sigma$; e.g. $\overline{Q}=(q_{\sigma(i)\sigma(j)})$.
Finally, denote by $\sigma_{a:b}$ any sub-vector $(\sigma(a),\ldots,\sigma(b)$)
of $\sigma$, where $a$, $b$ are positive integers. 

We show now that $(z_{t},m_{t})$ verify (\ref{eq:sys2}). Clearly,
$(z_{0},m_{0})=(1,1)$ with probability $1$. Then, compare\[
p(\sigma_{2:K}=\tau|z_{1}=2,s_{0}=k,\theta)=p(\sigma_{2:K}=\tau|s_{1}\neq k,s_{0}=k,\theta)\]
with \[
p(\sigma_{2:K}=\tau|z_{1}=1,s_{0}=k,\theta)=p(\sigma_{2:K}=\tau|s_{0}=s_{1}=k,\theta),\]
where $\tau$ is a vector of size $K-1$, and '1' and `2' are the
two only possible values for $z_{1}$. Since $(s_{t})$ is Markov,
conditional on $s_{0}$ the order of appearance $\sigma_{2:K}$ of
the states that differ from $s_{0}$ does not depend on the (random)
date $\eta$ where $s_{t}$ changes value for the first time; i.e.
$s_{0}=s_{1}=\ldots=s_{\eta-1}\neq s_{\eta}$. In particular, both
probabilities above equal \[
p(\sigma_{2:K}=\tau|s_{0}=k,\theta)\]
and one deduces that, conditional on $\sigma(1)$ (which equals $s_{0}$)
and $\theta$, $\sigma_{2:K}$ and $z_{1}$ are independent random
variables. Thus\[
p(z_{1}=1|\sigma,\theta)=p(s_{1}=s_{0}|s_{0},\theta)=q_{s_{0}s_{0}}=\overline{q}_{11}\]
almost surely. Since $\overline{\theta}$ is a deterministic function
of $\sigma$ and $\theta$, and $\overline{q}_{11}$ is a deterministic
function of $\theta$, one has: \[
p(z_{1}=1|\overline{\theta})=\overline{q}_{11},\]
and, \[
p(z_{1}=2|\overline{\theta})=1-p(z_{1}=1|\overline{\theta})=\sum_{i=2}^{K}\overline{q}_{1i}.\]
 This reasoning can be generalised to further time steps: for $k,l\leq m$,
where $m$ is the length of integer sequence $s$, $\sigma_{m+1:K}$
and $z_{t+1}$ are independent variables, conditional on $z_{t}=k$,
$m_{t}=m$, and $\sigma_{1:m}$. Thus, for $k,l\leq m$, and for some
arbitrary permutation $\tau$, \begin{eqnarray}
\lefteqn{p(z_{t+1}=l|z_{t}=k,m_{t}=m,\sigma=\tau,\theta)}\label{eq:toto}\\
 & = & p(s_{t+1}=\tau(l)|s_{t}=\tau(k),m_{t}=m,\sigma=\tau,\theta)\nonumber \\
 & = & p(s_{t+1}=\tau(l)|s_{t}=\tau(k),m_{t}=m,\sigma_{1:m}=\tau_{1:m},\theta)\nonumber \\
 & = & q_{\tau(k)\tau(l)}\nonumber \end{eqnarray}
which gives: \[
p(z_{t+1}=l|z_{t}=k,m_{t}=m,\sigma,\theta)=q_{\sigma(k)\sigma(l)}=\overline{q}_{kl},\]

with probability one, and, since probabilities sum up to one, \begin{equation}
p(z_{t+1}=m+1|z_{t}=k,m_{t}=m,\sigma,\theta)=\sum_{l=m+1}^{K}\overline{q}_{kl}.\label{eq:tata}\end{equation}
One can replace $(\sigma,\theta)$ with $\overline{\theta}$ in the
conditioning of (\ref{eq:toto}) and (\ref{eq:tata}), since the right-hand
sides depend only on $\overline{\theta}$, a deterministic function
of $(\sigma,\theta)$. This gives the desired result.

\section{Frequentist equivalence}

The two formulations define two incompatible data generating processes.
For instance, under the first model, the marginal distribution of
$y_{0}$ is a mixture:\[
p(y_{0}|\theta)=\sum_{k=1}^{K}\vartheta_{k}f(y_{0}|\xi_{k}),\]
whereas, under the second model, this distribution is $f(y_{0}|\overline{\xi}_{1})$.
However, this distinction seems relevant only when one would observe
repeated samples $y_{1:T}^{(i)}$ for the same parameter value $\theta$,
which seldom happens in practice. More importantly, the two formulations
define non equivalent likelihood functions, in the sense that, for
$t\in\Theta$, the likelihood of the first model $p(y_{1:T}|\theta=t)$
is generally not equal to the likelihood of the second model $p(y_{1:T}|\overline{\theta}=t)$,
if only because the former likelihood is invariant with respect to
state relabelling, whereas the latter is not. 

In fact, our time-ordered model corresponds to a re-parametrisation
of the full vector of unknowns, that is, $(\theta,s_{1:T})$ transformed
into $(\overline{\theta},z_{1:T})$, and therefore the complete likelihood
$p(y_{1:T}z_{1:t}|\overline{\theta})$ is a re-parametrised version
of $p(y_{1:T}s_{1:T}|\theta)$, but this property does not extend
to the marginal likelihood function $p(y_{1:T}|\theta)$. A practical
consequence is that maximum likelihood (or similar) estimators obtained
from either formulations cannot be compared easily; thus one should
refrain in principle from applying Frequentist estimation procedures
to sequentially ordered hidden Markov models.

\section{Discussion and extensions }

This example shows that Bayesian re-parametrisation is more powerful
a concept than Frequentist re-parametrisation. Since Bayesian analysis
treats equally all unknown quantities as random variables, whether
they are `parameters' or `latent variables', one can re-parametrise
(apply a one-to-one transform to) the full vector of unknowns, i.e.
$(\theta,z)$ above, rather than re-parametrise only the parameter
vector $\theta$. This explains why both formulations are equivalent
in terms of posterior distributions, but not in terms of (marginal)
likelihood functions. The author believes that this is yet another
example of the greater internal consistency of the Bayesian approach. 

One can easily derive sequentially ordered formulations for models
closely related to hidden Markov models, e.g. hidden semi-Markov models,
or continuous-time jump Markov models. The proof of the validity of
these sequentially ordered formulations is a simple extension of Section
1. For instance, a hidden Markov model is a semi-Markov model where
times between changes are geometrically distributed, but since our
proof is not based on this particular assumption, it extends readily
to semi-Markov models, up to cosmetic changes in the notations. In
the same way, a continuous-time jump Markov process is a process that
stays in a given state $k$, for a random, exponentially-distributed
duration, and then switch to another state according to some probability
transition matrix with zeros on its main diagonal. 

\bibliographystyle{apalike}
\bibliography{complete}

\end{document}